\input amstex
\documentstyle{amsppt}
\magnification=\magstep1
\pageheight{23truecm}\pagewidth{16truecm}
\vcorrection{0.5cm}
\parindent=10pt
\TagsOnRight
\document
\topmatter
\title
Automorphisms of free product-type \\
and their crossed-products
\endtitle
\author
Fumio HIAI$^{*}$ and Yoshimichi UEDA$^{**}$
\endauthor
\affil
${}^{*}$Graduate School of Information Sciences \\
Tohoku University \\
Aoba-ku, Sendai 980-8579, Japan \\
e-mail: hiai\@math.is.tohoku.ac.jp \\
\\
$^{**}$
Department of Mathematics \\
Graduate School of Science \\
Hiroshima University \\
Higashi-Hiroshima, 739-8526, Japan \\
e-mail: ueda\@math.sci.hiroshima-u.ac.jp
\endaffil
\thanks
Y.U.'s current address: Faculty of Mathematics, Kyushu University, 
Fukuoka, 810-8560, Japan. \newline e-mail: ueda\@math.kyushu-u.ac.jp
\endthanks 
\abstract
A continuous family of non-outer conjugate aperiodic automorphisms
whose crossed-products are all isomorphic is given on every
interpolated free group factor. An explicit ``duality" relationship
between compact co-commutative Kac algebra (minimal) free product
actions and free shift actions is also discussed.
\endabstract
\endtopmatter
\NoRunningHeads
\document
\subhead
1. Introduction
\endsubhead

\medskip
In this notes, we will study, as a natural continuation of [U1],
actions of the following type: Let $P$, $Q$ and $N$ be von Neumann
algebras and $\alpha$ be an action on $N$ of a group-like object $G$
such as groups or Kac algebras. Our objective here is to analyze the
natural extended action $\widetilde{\alpha}$ of $G$ on the free
product von  Neumann algebra $M := (P\otimes N) * Q$ defined by
$\widetilde{\alpha} = (\text{Id}_P\otimes\alpha) * \text{Id}_Q$.

As a simple application, we would like to point out that every
interpolated free group factor $L({\Bbb F}_r)$ (with arbitrary $r \in
(1,\infty]$) has continuously many non-outer conjugate aperiodic
automorphisms whose crossed-products are all isomorphic. This is a simple
supplementary remark to a famous result of J.~Phillips [Ph1] (also see [Ph2]). 
The same is true for a large class of free Araki-Woods factors (including 
all the unique type III$_\lambda$ cases with $\lambda \neq 0,1$) too. 
We also discuss an explicit relationship between the (minimal) action of 
compact co-commutative Kac algebra ${\Bbb K}_G$ (coming from a discrete group 
$G$) considered in [U1] and a certain $G$-free shift action. This simple but 
interesting observation by the first-named author was the starting point 
of this joint work. 

Although this work was carried out some time ago, typing of the notes
was done in part while the second-named author's stay in MSRI.
Y\.U\. would like to thank MSRI for their hospitality and financial
support.

\bigskip
\subhead
2. Aperiodic automorphisms of interpolated free group factors
\endsubhead

\medskip
Let us assume that $(Q, \tau_Q)$ is one of the following:
\roster
\item an interpolated free group factor $L({\Bbb F}_r)$ ($1 < r \leq
        \infty$) with the unique tracial state,
\item the group von Neumann algebra $L({\Bbb Z})$ with the canonical
        tracial state,
\item a finite direct sum of the trivial algebra ${\bold C}$ with
        a faithful (tracial) state, whose free dimension (see [D1], [D3])
        is in $(0,1)$.
\endroster
Then we will consider the free product (see [VDN])
$$
(M, \tau_M) := (R\otimes L^{\infty}({\Bbb T}), \tau_R\otimes\mu)
* (Q, \tau_Q), \tag{2.1}
$$
where $R$ is the hyperfinite type II$_1$ factor and $\mu$ is 
the expectation by the Haar probability measure on the 1-dimensional 
torus ${\Bbb T}$. A result of K.~Dykema (see [D1]) says that every 
interpolated free group factor (see [D2], [R1]) can be realized 
in this way. Note that the  algebra $L^{\infty}({\Bbb T})$ has
the canonical generator $u(z) = z$ ($z \in {\Bbb T}$), being a Haar
unitary, i.e., $\mu(u^n) = 0$ as long as $n \neq 0$.

Then the free product von Neumann algebra $M$ has the following natural
action of the 1-dimensional torus ${\Bbb T}$:
$$
\alpha_z |_{R \cup Q} := \text{the trivial action}, \quad
\alpha_z(u) := z u\ \text{(the multiplication of $z$)}
$$
for $z \in {\Bbb T}$. This action is nothing less than a free product
action considered in [U1], [U5] (also see [SU]) for the other purpose.

\proclaim{Lemma 1}
{\rm (1)} The action $\alpha$ is continuous in $p$-topology
{\rm (}or equivalently, in $u$-topology{\rm )}.

\noindent
{\rm (2)} The action $\alpha$ is outer, that is,
$\alpha_z \notin \text{\rm Int}(M)$ for any $z \neq 1$.
\endproclaim
\demo{Proof} (1) Straightforward.

(2) It suffices to show that $(M^{(\alpha,{\Bbb T})})'\cap M =
{\bold C}1$, where $M^{(\alpha,{\Bbb T})}$ denotes the fixed-point
algebra of the action $\alpha$ of $\Bbb T$. By [U3, Corollary 2]
(based on the idea of [P1, Lemma 2.5], or the use of [P2, Theorem 4.1]),
we have
$(R\otimes{\bold C}1)' \cap M = {\bold C}1_R\otimes L^{\infty}({\Bbb T})$.
Since both the free components do never commute with and since
$R\otimes{\bold C}$ and $Q$ sit in $M^{(\alpha,{\Bbb T})}$, the
relative commutant in question must be trivial.
\qed
\enddemo

\noindent
{\bf Remark 2.} (The free group factor version of Blattner's result [B])
\enspace Let $G$ be a separable locally compact group $G$. If the torus
action on $L^{\infty}({\Bbb T})$ was replaced, in the above construction,
by a faithful action of $G$ on a hyperfinite  von Neumann algebra, then
one would get an outer (continuous) action of $G$ on an interpolated free
group factor $L({\Bbb F}_r)$ with arbitrary $r \in (0,\infty]$. Note here
that R.-J. Blattner [B] (also see [T, p.\,47]) showed that any
separable locally compact group can outerly act on the hyperfinite type
II$_1$ factor $R$ so that the construction does really work. Hence,
Blattner's result mentioned just above remains still valid for any
interpolated free group factor. This observation was the initial
motivation of the work [U1].

\proclaim{Lemma 3} {\rm (cf\. [Ph1])}\enspace
Let $z_1, z_2 \in {\Bbb T}$ be irrationals, i.e.,
$z_k = e^{2\pi i \theta_k}$ with $\theta_k \notin {\Bbb Q}$
{\rm (}$k = 1,2${\rm )}. If neither $z_1=z_2$ nor $z_1=\bar z_2$
holds, then the corresponding aperiodic automorphisms
$\alpha_{z_1}$ and $\alpha_{z_2}$ never become outer conjugate to each
other.
\endproclaim
\demo{Proof} We make use of the idea given in [Ph1].
Namely the $\tau$-invariant (an invariant for outer conjugacy introduced
in [C1] for modular actions) will be used to distinguish.
The $\tau$-invariant $\tau(M,\alpha_{z_{k}})$ ($k = 1,2$) is the weakest
topology on ${\Bbb Z}$ making the mapping $n \mapsto \alpha_{z_k}^n \in
\text{Out}(M)$ continuous, where the quotient
$\text{Out}(M)=\text{Aut}(M)/\text{Int}(M)$ has the perfect
topological sense since $M$ is known to be a full factor. It is
straightforward to see that $\tau(M,\alpha_{z_k})$ is captured as  the
weak topology on
${\Bbb Z}$ induced by the mapping $n \mapsto z_k^n \in {\Bbb T}$ thanks
to Lemma 1\, (1). Hence [Ph1, Lemma 1.3] does work. Here, note
that $z_1\ne z_2,\bar z_2$ means
$\{z_1^{n} : n \in {\Bbb Z}\} \neq \{z_2^n : n \in {\Bbb Z}\}$.
\qed
\enddemo

Here we should mention that the $\tau$-invariant has been examined by 
D.~Shlyakhtenko [S3,\S8] for some automorphisms on free products of 
von Neumann algebras. The result there may be useful for further 
investigations. 

\proclaim{Lemma 4} {\rm (cf\. [U4, Proposition 1])}\enspace
If $z$ is irrational, i.e., $z = e^{2\pi i \theta}$  with $\theta
\notin {\Bbb Q}$, then the crossed-product $M\rtimes_{\alpha_z}{\Bbb Z}$
is isomorphic to the amalgamated free product {\rm(}see
{\rm [P2], [U2], [VDN])}
$$
(R\otimes(L^{\infty}({\Bbb T})\rtimes_{\alpha_z}{\Bbb Z})) *_{{\bold
C}\rtimes{\Bbb Z}}(Q\otimes L({\Bbb Z})), \tag{2.2}
$$
and its isomorphism class does not depend on the choice of $z$
{\rm (}or say $\theta${\rm )}. Here, the amalgamated free product in
\thetag{2.2} is taken with respect to the conditional expectations
$\tau_R\otimes E_{\theta}$ and $\tau_Q\otimes \text{\rm Id}_{L({\Bbb Z})}$
with the canonical conditional expectation
$E_{\theta} : L^{\infty}({\Bbb T})\rtimes_{\alpha_z}{\Bbb Z}
\rightarrow {\bold C}\rtimes{\Bbb Z}$.
\endproclaim
\demo{Proof} By [U4, Proposition 1], it suffices to show that
the amalgamated free product in \thetag{2.2} does not depend on
the choice of $z$. We first note that
$$
L^{\infty}({\Bbb T})\rtimes_{\alpha_z}{\Bbb Z} =
\{ u, v\}'' \cong  A_{\theta}'',
$$
where $v$ denotes the generator of ${\Bbb Z}$ in the crossed-product
so that $vu = z uv = e^{2\pi i \theta} uv$,
and $A_{\theta}''$ means the weak-closure of the irrational rotation
algebra $A_{\theta}$ via the GNS-representation associated with
the unique tracial state. Both the abelian subalgebras
$\{u\}''$ and $\{v\}''$ are known to be Cartan subalgebras in
$A_{\theta}''$, and hence the inclusion
$$
L^{\infty}({\Bbb T})\rtimes_{\alpha_z}{\Bbb Z} \supseteq
{\bold C}\rtimes{\Bbb Z}
$$
forms a pair of the hyperfinite type II$_1$ factor and a Cartan
subalgebra.
A.~Connes, J.~Feldman and B.~Weiss' result [CFW] says that all the 
Cartan subalgebras in any fixed hyperfinite factor are conjugate, 
which implies the assertion. \qed
\enddemo

Summing up the above three lemmas, we get the following theorem.

\proclaim{Theorem 5} Every interpolated free group factor has
continuously many non-outer conjugate aperiodic automorphisms whose
crossed-products are all isomorphic type II$_1$ factors.
\endproclaim

\medskip
We further investigate the crossed-product of $M$ by an
aperiodic automorphism of the form $\alpha_z$ with irrational
$z \in {\Bbb T}$. Thanks to Lemma 4 (and the proof), we may
investigate the following amalgamated free product:
$$
N := (R\otimes R_0) *_D (Q\otimes D),
$$
where $R_0 \supseteq D$ is a (unique) pair of the hyperfinite type
II$_1$ factor and a Cartan subalgebra.

\proclaim{Theorem 6} {\rm (cf\. [U5, Theorem 8])}\enspace
Let $\omega$ be a free ultrafilter. Then the crossed-product
$N = M\rtimes_{\alpha_z}{\Bbb Z}$ with irrational $z$ satisfies
$$
N_{\omega} = N'\cap N^{\omega} = R_0{}' \cap D^{\omega}, \tag{2.3}
$$
which implies that $N$ has the property $\Gamma$ {\rm(}i.e.,
not full{\rm\,)} but not McDuff {\rm ([Mc])}.
\endproclaim
\demo{Proof} We make use of the same idea as in [U5]. Since $R$ is
a type II$_1$ factor, we can choose a unitary $u$ in $R$
($=R\otimes{\bold C}1\subseteq M$) with
$\tau_R(u^n) = 0$ as long as
$n \neq 0$. Towards making use of the above-mentioned idea, we have to
choose an invertible element $y$ in $Q$ with $\tau_Q(y) = 0$. If
$(Q,\tau_Q)$ is as in the case \thetag{1} or \thetag{2},
then one can choose a Haar unitary so that there is no problem.
When $(Q,\tau_Q)$ is as in the case \thetag{3}, we write
$$
Q =
{\bold C}\oplus{\bold C}\oplus\cdots\oplus{\bold C} \quad
(\text{$n$ times}),
$$
and the (faithful) tracial state comes from a vector of weights
$\lambda = (\lambda_1,\lambda_2,\dots,\lambda_n)$ with $\lambda_i > 0$,
$\sum_{i = 1}^n \lambda_i = 1$. The above invertible $y$ should be
a vector $\mu = (\mu_1,\mu_2,\dots,\mu_n)$ in ${\bold C}^n \cong Q$
such that all the $\mu_i$'s are non-zero and $\mu$ is
orthogonal to $\lambda$ with respect to the usual inner product on
${\bold C}^n$.  For choosing such a vector $\mu$, it suffices to
show that the complement $({\bold C}\lambda)^{\perp}$ is not contained in
any proper subspace generated by a part of the standard basis
$\{e_1,e_2,\dots, e_n\}$ of ${\bold C}^n$. If so was, then there would be
a subset $\{e_{i_1},\dots,e_{i_j}\}$ such that $[e_{i_1},\dots,e_{i_j}]
\subseteq (({\bold C}\lambda)^{\perp})^{\perp} = {\bold C}\lambda$,
which is a contradiction to the condition that all the $\lambda_i \neq 0$.
Hence we have shown the existence of an invertible $y$ in $Q$ with
$\tau_Q(y) = 0$.

As in [U5, Proposition 5], we have, for each
$x \in \{u\}'\cap N^{\omega}$,
$$
\|y(x - (E_D^N)^{\omega}(x))\|_{\tau_N^{\omega}} \leq
\|yx - xy\|_{\tau_N^{\omega}},
$$
and hence if $x$ furthermore commutes with $y$, then
$x = (E_D^N)^{\omega}(x) \in D^{\omega}$ since $y$ is invertible.
Therefore, we have obtained the formula \thetag{2.3}.

The relative commutant $R_0{}'\cap D^{\omega}$ can be identified
(abstractly) with the von Neumann algebra generated by equivalence classes
of $\omega$-centralizing sequences (under the action of an ergodic
finite-measure preserving transformation) consisting of Borel subsets
in a non-atomic Lebesgue space $X$ with $D = L^{\infty}(X)$. In this case,
the set of those $\omega$-centralizing sequences is known to be very large
(see [Sc, Theorems 3.1, 3.3]), and thus the latter assertion follows. \qed
\enddemo

\medskip
The pair $R_0 \supseteq D$ is known to be constructed from
the type II$_1$ amenable discrete equivalence relation ${\Cal R}$
over a non-atomic Lebesgue space $X$, that is,
$$
R_0 = W^*({\Cal R}) \supseteq D = L^{\infty}(X).
$$
(See [CFW], [FM].) The Galois group $\text{Gal}(R_0\supseteq D) :=
\{\alpha \in \text{Aut}(R_0) : \alpha|_D = \text{Id}\}$ is identified
with the group $Z^1({\Cal R};{\Bbb T})$ of 1-cocycles on ${\Cal R}$
via the mapping
$c\in Z^1({\Cal R};{\Bbb T}) \mapsto M_c\in \text{Gal}(R_0 \supseteq D)$,
where $M_c$ is the multiplier defined on $R_0 = W^*({\Cal R})$ (see [FM]).
We define the homomorphism $\Phi : Z^1({\Cal R};{\Bbb T}) \rightarrow
\text{Aut}(N)$ (or more precisely into $\text{Gal}(N \supseteq D)$)
by
$$
\Phi(c) := (\text{Id}_R\otimes M_c) *_D (\text{Id}_Q\otimes\text{Id}_D),
\quad c \in Z^1({\Cal R};{\Bbb T}).
$$
We first point out that $\Phi(B^1({\Cal R};{\Bbb T})) = \text{Int}(N,D)$
($:= \{\text{Ad}u \in \text{Aut}(N) : u \in {\Cal U}(D)\}$), and that
if $\beta$ is an approximate inner automorphism on $N$, then it must be
of the form: $\alpha = \text{Ad}X\circ\beta_0$ with $X \in {\Cal U}(N)$,
$\beta_0 \in \overline{\text{Int}(N,D)}$ thanks to Theorem 6 and
Connes' method (see [U5, \S\S2.3] for summary). Then $\beta_0$
should act on both $R$ and $Q$ trivially and satisfy that its restriction
to $R_0$ is in $\overline{\text{Int}(R_0, D)}$, and hence
$\beta_0 = (\text{Id}_R\otimes M_c) *_D (\text{Id}_Q\otimes\text{Id}_D)
= \Phi(c)$
for some $c \in \overline{B^1({\Cal R}; {\Bbb T})}$. What we explained
here is exactly the same as in [U5, \S4], and we get the
following  proposition in this way, see [U5, \S4] for details.

\proclaim{Proposition 7} {\rm (cf\. [U5, Theorem 14])}\enspace
The crossed-product $N = M\rtimes_{\alpha_z}{\Bbb Z}$ with irrational $z$
satisfies
$$
\chi(N) = \overline{\text{\rm Int}(N)}/\text{\rm Int}(N) \cong
H^1({\Cal R}; {\Bbb T}),
$$
and the isomorphism is induced by $\Phi$. Here, $\chi(N)$ is the
$\chi$-group {\rm [C2]} and $H^1({\Cal R}; {\Bbb T})$ is the first
cohomology group of the type II$_1$ ergodic amenable discrete
equivalence relation {\rm (}see {\rm [Sc])}.
\endproclaim

The free Araki-Woods factor $\Gamma({\Cal H}_{\bold R}, U_t)''$ (equipped
with the so-called free quasi-free state $\varphi_U$) associated with
a (non-trivial) one-parameter group of orthogonal transformations on a
real Hilbert  space ${\Cal H}_{\bold R}$ was introduced by D.~Shlyakhtenko 
[S1]. It was shown in [S1], [S2] that $\Gamma({\Cal H}_{\bold R},U_t)''$ 
is a factor of type III$_\lambda$ ($0<\lambda<1$) if $U_t$ is periodic with
period $2\pi/\log\lambda$ and a factor of type III$_1$ if $U_t$ is
non-periodic. Replacing
$(Q,\tau_Q)$ in
\thetag{2.1} by
$(\Gamma({\Cal H}_{\bold R}, U_t)'',\varphi_U)$, we consider
$$
(M, \varphi) := (R\otimes L^{\infty}({\Bbb T}), \tau_R\otimes\mu)
* (\Gamma({\Cal H}_{\bold R}, U_t)'',\varphi_U)
$$
and aperiodic automorphisms $\alpha_z$ for irrational $z\in{\Bbb T}$
defined as above. If $U_t$ has an eigenvalue not equal to $1$, then
we notice, by [S1, Corollary 5.5] and [D1], the following
state-preserving isomorphisms:
$$
\align
(M,\varphi) &\cong
(R\otimes L^{\infty}({\Bbb T}), \tau_R\otimes \mu)
* (L({\Bbb F}_\infty),\tau) * (\Gamma({\Cal H}_{\bold R}, U_t)'',\varphi_U) \\
&\cong
(L({\Bbb F}_\infty),\tau) * (\Gamma({\Cal H}_{\bold R}, U_t)'',\varphi_U) \\
&\cong(\Gamma({\Cal H}_{\bold R}, U_t)'',\varphi_U),
\endalign
$$
and all the arguments above can work in this setting as well. With
$R_0 \supseteq D$ as above, we hence have
$$
M\rtimes_{\alpha_z}{\Bbb Z}\cong
(R\otimes R_0)*_D(\Gamma({\Cal H}_{\bold R}, U_t)''\otimes D),
$$
which is a factor of the same III$_{\lambda}$-type ($0<\lambda\le1$) as
$\Gamma({\Cal H}_{\bold R},U_t)''$ (see e.g.\ [U2, Theorem 2.6, Corollary
4.5]). The following theorem is obtained in this way.

\proclaim{Theorem 8}
The free Araki-Woods factor $\Gamma({\Cal H}_{\bold R}, U_t)''$ with $U_t$ 
having an eigenvalue not equal to $1$ has continuously many non-outer 
conjugate aperiodic automorphisms, all of whose crossed-products are 
isomorphic to a non-full factor of the same III$_{\lambda}$-type as 
$\Gamma({\Cal H}_{\bold R}, U_t)''$, not McDuff and having
the same $\chi$-group as in Proposition {\rm 7}.
\endproclaim

\bigskip
\subhead
3. Duality between free product actions and free shift actions
\endsubhead

\medskip
We are now in turn going to deal with actions whose crossed-products still
stay in the ``category'' consisting of interpolated free group factors.
It is probably well-known that free shift actions are such typical examples.

We would like here to point out an explicit relationship between the
action of a compact (co-commutative) Kac algebra ${\Bbb K}_G$ considered
in [U1] and a certain free shift action (or free permutation action
associated with $G$).

\medskip
Let $P$, $Q$ be von Neumann algebras with specific faithful normal
states $\varphi_P$, $\varphi_Q$, respectively. Let $G$ be a discrete
(countable) group, and $\lambda_g$ ($g \in G$) means its left regular
representation. Let $L(G)$ be the group von Neumann algebra with the
canonical trace $\tau_G$. The compact co-commutative Kac algebra
${\Bbb K}_G = (L(G),\delta_G,\kappa_G,\tau_G)$ can act on the free
product
$$
(M, \varphi) := (P\otimes L(G), \varphi_P\otimes\tau_G) * (Q,\varphi_Q)
$$
by the free product (co-)action
$$
\Gamma^G =
(\text{Id}_P\otimes\delta_G) * (\text{Id}_Q\otimes 1_{L(G)})
$$
in the sense of [U1]. This action is nothing less than that considered in
[U1, \S4].

Consider $P$, $Q$ and $L(G)$ as subalgebras of $M$ naturally, and set
$$
N := \Biggl( P \cup \bigcup_{g \in G} \lambda_g Q \lambda_g^* \Biggr)''.
$$
Since $\text{Ad}\,\lambda_g(x) = x$ ($x \in P$) and
$\text{Ad}\,\lambda_g(\lambda_h y\lambda_h^*) = \lambda_{gh} y
\lambda_{gh}^*$  ($y \in Q$, $h \in G$), we can define the action $\alpha$
of $G$ on $N$  by $\alpha_g := \text{Ad}\,\lambda_g|_N$ ($g \in G$).

It is rather easy to show that $\{P, \{\lambda_g Q \lambda_g^*\}_{g\in G}\}$
forms a free family with respect to $\varphi$ and that each $\lambda_g$
($g \neq e$) is orthogonal to $L^2(N,\varphi)$ in $L^2(M,\varphi)$. (See
[CD] for example, where a slightly generalized situation was treated in the
$C^*$-algebraic setting.) Note here that $\sigma_t^{\varphi}(\lambda_g) =
\lambda_g$ and $\sigma_t^{\varphi}(N) = N$ (globally invariant) so that
there exists the $\varphi$-conditional expectation $E:M\to N$ satisfying
$E(\lambda_g)=0$ ($g\ne e$). Therefore, we see that
$$
(N, \varphi|_N) \cong (P, \varphi_P) *
\left( \underset{g \in G}\to{*} (Q, \varphi_Q)_g \right),
\quad M \cong N\rtimes_{\alpha}G,
$$
and that $\alpha$ can be naturally identified with the free product
action $\text{Id}_P * \gamma^G$ under the above identification, where
$\gamma^G$ denotes the $G$-free shift action on
$\displaystyle{\underset{g \in G}\to{*}(Q, \varphi_Q)_g}$. What we want
to point out here is the following ``duality'' between the $G$-free
shift action $\gamma^G$ and the free product action $\Gamma^G$:

\proclaim{Theorem 9} In the current setting, the following
assertions hold\,{\rm :}
\roster
\item $(M,\Gamma^G) \cong (N\rtimes_{\alpha}G,\widehat{\alpha})$ with the
        dual action $\widehat{\alpha}$
        {\rm ($= \widehat{\text{Id}_P * \gamma^G}$)}
        {\rm (}see {\rm [NT])},
\item $N = M^{\Gamma^G}$,
\item $M\rtimes_{\Gamma^G}{\Bbb K}_G \cong N\otimes B(\ell^2(G))$.
\endroster
\endproclaim
\demo{Proof} (1) Under the identification $M = N\rtimes_{\alpha}G$
explained above, the dual action $\widehat{\alpha}$ acts on $M$ in
the following manner:
$$
\align
\widehat{\alpha}(y) &= y \otimes 1_{L(G)}, \quad y \in N, \\
\widehat{\alpha}(\lambda_g) &= \lambda_g \otimes \lambda_g, \quad g \in G.
\endalign
$$
On the other hand, the definition of $\delta_G$ gives
$$
\align
\Gamma^G(\lambda_g) &= \lambda_g\otimes\lambda_g, \quad g \in G, \\
\Gamma^G(x) &= x\otimes 1_{L(G)}, \quad x \in P, \\
\Gamma^G(\lambda_g y \lambda_g^*) &=
(\lambda_g y \lambda_g^*)\otimes 1_{L(G)}, \quad y \in Q,\ g \in G,
\endalign
$$
and hence $\Gamma^G = \widehat{\alpha}$ follows.

(2) comes from (1) and the well-known formula
$N = (N\rtimes_{\alpha}G)^{\widehat{\alpha}}$ (see [NT]).

(3) comes from (1), (2) and the Takesaki duality (see [NT]).
\qed
\enddemo

Furthermore, similarly to Lemma 1\,(2), one can see that if the
centralizer $P_{\varphi_P}$ is diffuse and $Q\ne{\bold C}$, then the
action $\alpha$ is outer.

\medskip
\noindent{\bf Remarks 10.}
(1) According to [U4, Proposition 1] (or the proof of
[SU, Proposition 3.2]) one has
$$
\align
M\rtimes_{\Gamma^G}{\Bbb K}_G &\cong
\left( (P\otimes L(G))\rtimes_{\text{Id}_P\otimes\delta_G}{\Bbb K}_G \right)
*_{\ell^{\infty}(G)}
\left( Q\rtimes_{\text{Id}_Q\otimes 1_{L(G)}}{\Bbb K}_G \right) \\
&\cong
\left(P\otimes B(\ell^2(G))\right) *_{\ell^{\infty}(G)}
\left(Q\otimes\ell^{\infty}(G)\right).
\endalign
$$
In the above, we omitted the mention of the specific conditional
expectations since they are clear from the context. (And, in what
follows, we will do omit it as long as when no confusion is possible.)
The above computation says the following (probably known) simple fact:
If $A$ is a (unique) atomic MASA in $B({\Cal H})$ with arbitrary
$\dim{\Cal H}$, then
$$
\left(P\otimes B({\Cal H})\right) *_A \left(Q\otimes A\right)
\cong
\left(P * \left( Q^{*(\text{dim}{\Cal H})} \right)\right)\otimes
B({\Cal H}).
$$
(Compare with [PS, Theorem 3.3].)
Note here that this can be directly shown by simple algebraic method
(without any random matrix type technique). What we explained here was
one of the motivations of the works [U4], [SU].

(2) Theorem 9 suggests us what ``free shift actions'' associated with
the duals of compact quantum groups (e.g., $\widehat{\text{\rm SU}_q(N)}$) 
should be. In fact, the theorem says that the minimal free product action 
of $\text{\rm SU}_q(2)$ given in [U4] ($P = {\bold C}$ in that case) can be 
regarded as the dual (co-)action of the ``free 
$\widehat{\text{\rm SU}_q(2)}$-shift'' action. We will return to this point 
of view in future.

\bigskip
As mentioned before, the action $\alpha$ is the free product action
$\text{Id}_P * \gamma^G$, and one has
$$
\align
M = (P\otimes L(G)) * Q
&\cong
\left( P * \left(\underset{g \in G}\to{*} Q_g \right)\right)
\rtimes_{\alpha}G \\
&\cong
\left(P\otimes L(G)\right) *_{L(G)}
\left(\left(\underset{g \in G}\to{*}Q_g\right)\rtimes_{\gamma^G}G\right)
\endalign
$$
with $Q_g := Q$ ($g \in G$). This in particular says that
$\displaystyle{\left(\underset{g \in G}\to{*}Q_g\right)\rtimes_{\gamma^G}
G\cong L(G) * Q}$, and hence we get
$$
\left(P\otimes L(G)\right) * Q \cong
\left(P\otimes L(G)\right) *_{L(G)} \left(L(G) * Q\right). \tag{3.1}
$$
In fact, we have the following slightly more general fact.

\proclaim{Fact 11} {\rm (A special case of [NSS, Proposition 2.7])}
\enspace Let $P$, $Q$ and $A$ be von Neumann algebras with faithful
normal states
$\varphi_P$, $\varphi_Q$ and $\varphi_A$, respectively. Then
$$
\left(P\otimes A, \varphi_P\otimes\varphi_A\right) *_A
\left(A * Q, \varphi_A * \varphi_Q\right) \cong
\left((P\otimes A) * Q, (\varphi_P\otimes\varphi_A) * \varphi_Q\right).
$$
\endproclaim

This can be shown by checking a suitable freeness; the proof is
straightforward. Quite recently, A.~Nica, D.~Shlyakhtenko and
R.~Speicher [NSS] showed a more general fact as a corollary of what
they developed on the operator-valued $R$-transform.

\medskip
For each amenable group $G$ and each $r\in(1,\infty]$, let us specify
$(P,\tau_P)$ and $(Q,\tau_Q)$ for which we have
$$
M \cong L({\Bbb F}_r), \quad N \cong L({\Bbb F}_{(r-1)|G|+1}).
$$
To do so, it suffices, for example, to assume that $P$ is the hyperfinite
type II$_1$ factor $R$, $(Q, \tau_Q)$ is one of those given at the
beginning of Section 2 and the free dimension of $(Q,\tau_Q)$ is $r-1$.
Then the above realization of $M,N$ can be easily seen by simple
computations of free dimensions based  on the results of K.~Dykema [D1].
In this way, we have obtained the following result.

\proclaim{Proposition 12} For each amenable
{\rm(}countable{\rm)} group
$G$ and each $r \in (0,\infty]$, there exists an outer action $\alpha$
of $G$ on $L({\Bbb F}_{(r-1)|G|+1})$ such that
$$
L({\Bbb F}_{(r-1)|G|+1})\rtimes_{\alpha}G \cong L({\Bbb F}_r).
$$
In particular, for each $r\in(1,\infty]$, there exists an aperiodic
automorphism $\alpha$ on $L({\Bbb F}_{\infty})$ such that
$$
L({\Bbb F}_{\infty})\rtimes_{\alpha}{\Bbb Z} \cong L({\Bbb F}_r).
$$
\endproclaim

\medskip\noindent
{\bf Remark 13.}
When $G$ is an amenable group and $P,Q$ are specified as above, the
fact \thetag{3.1} says that the formula of free dimension for a certain
class of amalgamated free products (including the multi-matrix algebra
situation) given in [D3] is valid even in a case where the amalgamation
subalgebra is non-atomic.

\medskip
Finally, we would like to give small comments on the free analogs of
Bernoulli shifts, introduced by replacing the tensor product ``$\otimes$"
by the free product ``$*$".

The ``free Bernoulli shift" consisting of the free group factor
$L({\Bbb F}_{\infty})$ and its ergodic, aperiodic automorphism $\sigma_p$
associated with a (non-degenerate) probability vector $p=(p_1,\dots,p_n)$
is the free analog of  Connes-St{\o}rmer's Bernoulli shift ([CS]), and it
is constructed in the following manner (see e.g.\ [P3]):  Let
$$
(N, \psi) := \underset{k\in{\Bbb Z}}\to{*}
(M_n({\bold C}), \varphi_p)_k \tag{3.2}
$$
be the free product, where the state $\varphi_p$ has the
diagonal density matrix $\text{diag}(p_1,\dots,p_n)$. Then the ${\Bbb Z}$-free 
shift $\alpha$ preserves the free product state $\psi$ so that one can consider 
its restriction $\sigma_p:=\alpha|_{N_\psi}$ to the centralizer
$N_{\psi}$. (Note that, in the case of equal probabilities, i.e.,
$p_1 =\cdots = p_n =\frac{1}{n}$, the state $\psi$ itself is a trace and
$N\cong L({\Bbb F}_{\infty})$ due to [D1] so that we simply set
$\sigma_p:=\alpha$.) According to [D4], the centralizer
$N_{\psi}$ is isomorphic to $L({\Bbb F}_{\infty})$ (when $\psi$ is not a
trace).  The discussion given before Theorem 9 tells us that the
crossed-product $N\rtimes_{\alpha}{\Bbb Z}$ is isomorphic to the free
product
$$
(M, \varphi) := (L({\Bbb Z}),\tau_{\Bbb Z}) * (M_n({\bold C}),\varphi_p),
$$
which is known to be isomorphic to $L({\Bbb F}_{2-\frac{1}{n^2}})$ in
the case of equal probabilities (see [D1]), or otherwise to the free
Araki-Woods factor whose Connes' $Sd$ invariant is the multiplicative
group generated by $p_i/p_j$ ($1\le i,j\le n$) (see [S1, p.\,365]).
In fact, $N$, $\psi$ and $\alpha$ are realized as
$$
N=\Biggl(\bigcup_{k\in{\Bbb Z}}\lambda^k M_n({\bold C})\lambda^{*k}
\Biggr)''\,\bigl(\subseteq M\bigr),\quad \psi=\varphi|_N,\quad
\alpha=\text{Ad}\,\lambda,
$$
where $\lambda$ is the generating unitary of $L({\Bbb Z})$. Moreover,
note that $\varphi=\psi\circ E$ and $\sigma_t^\varphi|_N=\sigma_t^\psi$,
where $E:M\,(\cong N\rtimes_\alpha{\Bbb Z})\to N$ is the canonical
conditional expectation. Now, let $E_\varphi$ be the $\varphi$-conditional
expectation from $M$ onto the centralizer $M_\varphi$, and $E_\psi$ the
$\psi$-conditional expectation from $N$ onto $N_\psi$. Then it is easy to
see that $E_\varphi(\lambda)=\lambda$ (thanks to $\lambda\in M_\varphi$)
and $E_\varphi|_N=E_\psi$ (by the commuting square property thanks to
$E(M_\varphi)\subseteq N_\psi$). Therefore, we notice
$$
M_\varphi=(N_\psi\cup\{\lambda\})''\cong N_\psi\rtimes_{\sigma_p}{\Bbb Z}.
$$
Since $M_\varphi\cong L({\Bbb F}_\infty)$ (see [D4], [R2], [S1], [S2]),
we have  obtained the following proposition.

\proclaim{Proposition 14} Let $(L({\Bbb F}_{\infty}),\sigma_p)$ be the
free Bernoulli shift associated with a probability vector
$p=(p_1,\dots,p_n)$ defined above. Then
$$
L({\Bbb F}_{\infty})\rtimes_{\sigma_p}{\Bbb Z} \cong
\cases
L({\Bbb F}_{2-\frac{1}{n^2}}) &
\text{if $p_1 = \cdots = p_n = \frac{1}{n}$,} \\
L({\Bbb F}_{\infty}) &
\text{otherwise.}
\endcases
$$
\endproclaim

One can also define a more direct free analog of the classical Bernoulli
shift, which may be called the free ``commutative" Bernoulli shift
associated with a probability vector $p=(p_1,\dots,p_n)$. It consists
of $L({\Bbb F}_{\infty})$ and its ergodic, aperiodic automorphism
$\gamma_p$ again, provided by replacing the full matrix algebra
$M_n({\bold C})$ in \thetag{3.2} by its diagonal subalgebra
${\bold C}^n$. Since no type III situation appears in this case, we
simply define $\gamma_p$ by the ${\Bbb Z}$-free shift $\alpha$.
Similarly to the above consideration, we see that the crossed-product
$L({\Bbb F}_{\infty})\rtimes_{\gamma_p}{\Bbb Z}$ is isomorphic to the
free product
$$
(L({\Bbb Z}),\tau_{\Bbb Z}) * ({\bold C}^n, p),
$$
which is isomorphic to $L({\Bbb F}_{2-\|p\|_2^2})$ due to [D1]. Here,
$\|p\|_2$ denotes the $\ell^2$-norm of the probability vector $p$.
Therefore, the following proposition is obtained.

\proclaim{Proposition 15} Let $(L({\Bbb F}_{\infty}),\gamma_p)$ be the
above free Bernoulli shift associated with a probability vector
$p = (p_1,\dots,p_n)$. Then
$$
L({\Bbb F}_{\infty})\rtimes_{\gamma_p}{\Bbb Z}
\cong L({\Bbb F}_{2-\|p\|_2^2}),
$$
where $\|p\|_2$ is the $\ell^2$-norm of $p$.
\endproclaim

\medskip\noindent
{\bf Remarks 16.}
(1) It is interesting to note that aperiodic automorphisms on
$L({\Bbb F}_{\infty})$ given in this section have the characteristic very
different from that of automorphisms given in Section 2. More precisely,
if $\alpha$ is a ${\Bbb Z}$-free shift on $L({\Bbb F}_\infty)$ as in
Propositions 12, 14 and 15, then the crossed-product
$L({\Bbb F}_\infty)\rtimes_\alpha{\Bbb Z}$ is a certain interpolated free
group factor and hence full. This means (see [J], [Ph1], [Ph2]) that
$\{\alpha^n:n\in{\Bbb Z}\}$ forms a discrete subgroup of 
$\text{Out}(L({\Bbb F}_\infty))$. On the other hand, aperiodic
automorphisms in Theorem 5 have non-full crossed-products as stated in
Theorem 6.

(2) As proved in [D2], [R1] independently, the interpolated free group
factors are either all isomorphic or all non-isomorphic (i.e.,
$L({\Bbb F}_r)\not\cong L({\Bbb F}_{r'})$ for $r\neq r'$). If the latter
case were true, then one would have a continuous family of non-cocycle
conjugate (or equivalently non-outer conjugate) ${\Bbb Z}$-free shifts on
$L({\Bbb F}_\infty)$ by Proposition 12 as well as some cocycle conjugacy
classification results for free Bernoulli shifts $\sigma_p$ and $\gamma_p$
given in Propositions 14 and 15. In this connection, it may be pointed out
that the classification up to conjugacy by means of Connes-St{\o}rmer's
dynamical entropy ([CS]) is meaningless for these free Bernoulli shifts
since all such free shifts have zero entropy (see e.g.\ [St]).

\bigskip
\subhead
References
\endsubhead

\medskip\noindent
[B] R.-J.\ Blattner, Automorphic group representations,
{\it Pacific J.\ Math.}\ {\bf 8} (1958) 665--677.

\medskip\noindent
[C1] A.\ Connes, Almost periodic states and factors of type III$_1$,
{\it J.\ Funct.\ Anal.}\ {\bf 16} (1974), 415--455.

\medskip\noindent
[C2] A.\ Connes, Sur la classification des factors de type II,
{\it C.\ R.\ Acad.\ Sci.\ Paris S\'er.\ A} {\bf 281} (1975), 13--15.

\medskip\noindent
[CD]  M.\ Choda and K.-J.\ Dykema, Purely infinite, simple $C^*$-algebras
arising from free product constructions.\ III,
{\it Proc.\ Amer.\ Math.\ Soc.}\ {\bf 128} (2000), 3269--3273.

\medskip\noindent
[CFW] A.\ Connes, J.\ Feldman and B.\ Weiss, An amenable equivalence
relation is generated by a single transformation,
{\it Ergodic Theory Dynamical Systems} {\bf 1} (1981), 431--450.

\medskip\noindent
[CS] A.\ Connes and E.\ St{\o}rmer, Entropy for automorphisms of II$_1$
von Neumann algebras,
{\it Acta Math.}\ {\bf 134} (1974), 289--306.

\medskip\noindent
[D1] K.\ Dykema, Free products of hyperfinite von Neumann algebras and
free dimension,
{\it Duke Math.\ J.}\ {\bf 69} (1993), 97--119.

\medskip\noindent
[D2] K.\ Dykema, Interpolated free group factors,
{\it Pacific J.\ Math.}\ {\bf 163} (1994), 123--135.

\medskip\noindent
[D3] K.\ Dykema, Amalgamated free products of multi-matrix algebras
and a construction of subfactors of a free group factor,
{\it Amer.\ J.\ Math.}\ {\bf 117} (1995), 1555--1602.

\medskip\noindent
[D4] K.\ Dykema, Free products of finite dimensional and other
von Neumann algebras with respect to non-tracial states,
in {\it Free probability theory},
Fields Inst.\ Commun., Vol.\ 12, Amer.\ Math.\ Soc.,
Providence, RI, 1997, pp.~41--88.

\medskip\noindent
[FM] J.\ Feldman and C.\ C.\ Moore, Ergodic equivalence relations,
cohomology, and von Neumann algebras.\ I, II,
{\it Trans.\ Amer.\ Math.\ Soc.}\ {\bf 234} (1977), 289--324, 325--359.

\medskip\noindent
[J] V.\ F.\ R.\ Jones, Central sequences in crossed products of
full factors, {\it Duke Math.\ J.}\ {\bf 49} (1982), 29--33.

\medskip\noindent
[Mc] D.\ McDuff, Central sequences and the hyperfinite factor,
{\it Proc.\ London Math.\ Soc.\ (3)} {\bf 21} (1970), 443--461.

\medskip\noindent
[NSS] A.\ Nica, D.\ Shlyakhtenko and R.\ Speicher,
A characterization of freeness by a factorization property of
$R$-transform, Preprint (2001).

\medskip\noindent
[NT] Y.\ Nakagami and M.\ Takesaki, Duality for Crossed Products of
von Neumann Algebras, Lecture Notes in Mathematics, Vol.\ 731, Springer,
Berlin, 1979.

\medskip\noindent
[Ph1] J.\ Phillips, Automorphisms of full II$_1$ factors, with
applications to factors of type III,
{\it Duke Math.\ J.}\ {\bf 43} (1976), 375--385.

\medskip\noindent
[Ph2] J.\ Phillips, Automorphisms of full II$_1$ factors, II, 
{\it Canad.\ Math.\ Bull.}, {\bf 21} (1978), 325--328.

\medskip\noindent
[P1] S.\ Popa, Orthogonal pairs of $*$-subalgebras in finite von Neumann
algebras, {\it J.\ Operator Theory} {\bf 9} (1983), 253--268.

\medskip\noindent
[P2] S.\ Popa, Markov traces on universal Jones algebras and subfactors
of finite index. {\it Invent.\ Math.}\ {\bf 111} (1993), 375--405.

\medskip\noindent
[P3] S.\ Popa, A rigidity result for actions of property T groups
by Bernoulli shifts, Preprint (2001).

\medskip\noindent
[PS] S.\ Popa and D.\ Shlyakhtenko, Universal properties of
$L({\Bbb F}_{\infty})$ in subfactor theory,
MSRI Preprint No.\ 2000-032.

\medskip\noindent
[R1] F.\ R\u adulescu, Random matrices, amalgamated free products and
subfactors of the von Neumann algebra of a free group, of noninteger
index, {\it Invent.\ Math.}\ {\bf 115} (1994), 347--389.

\medskip\noindent
[R2] F.\ R\u adulescu, A type $III_{\lambda}$ factor with core isomorphic
to the von Neumann algebra of a free group, tensor $B(H)$,
in {\it Recent Advances in Operator Algebras}, {\it Ast\'erisque} {\bf 232}
(1995), 203-209.

\medskip\noindent
[Sc] K.\ Schmidt, Algebraic Ideas in Ergodic Theory, CBMS Regional
Conference Series in Mathematics, Vol.\ 76, Amer. Math. Soc.,
Providence, RI, 1990.

\medskip\noindent
[S1] D.\ Shlyakhtenko, Free quasi-free states,
{\it Pacific J.\ Math.}\ {\bf 177} (1997), 329--368.

\medskip\noindent
[S2] D.\ Shlyakhtenko, Some applications of freeness with amalgamation,
{\it J.\ Reine Angew.\ Math.}\ {\bf 500} (1998), 191--212.

\medskip\noindent
[S3] D.\ Shlyakhtenko, {$A$}-valued semicircular systems, 
{\it J. Funct. Anal.}\ {\bf 166} (1999), 1--47

\medskip\noindent
[SU] D.\ Shlyakhtenko and Y.\ Ueda, Irreducible subfactors of
$L({\Bbb F}_{\infty})$ of index $\lambda > 4$, 
{\it J.\ Reine Angew.\ Math.} {\bf 548} (2002), 149--166. 

\medskip\noindent
[St] E.\ St{\o}rmer, States and shifts on infinite free products of
$C\sp *$-algebras, in {\it Free probability theory}, Fields Inst.\ Commun.,
Vol.\ 12, Amer.\ Math.\ Soc., Providence, RI, 1997, pp.~281--291.

\medskip\noindent
[T] M.\ Takesaki, Structure of Factors and Automorphism Groups, CBMS
Regional Conference Series in Mathematics, Vol.\ 51,
Amer.\ Math.\ Soc.\, Providence, RI, 1983.

\medskip\noindent
[U1] Y.\ Ueda, A minimal action of the compact quantum group 
$\text{\rm SU}_q(n)$ on a full factor,
{\it J.\ Math.\ Soc.\ Japan} {\bf 51} (1999), 449--461.

\medskip\noindent
[U2] Y.\ Ueda, Amalgamated free product over Cartan subalgebra,
{\it Pacific J.\ Math.}\ {\bf 191}, (1999), 359--392.

\medskip\noindent
[U3] Y.\ Ueda, Remarks on free products with respect to non-tracial
states, 
{\it Math.\ Scand.}\ {\bf 88}, (2001), No. 1, 111--125. 

\medskip\noindent
[U4] Y.\ Ueda, On the fixed-point algebra under a minimal free
product-type action of the quantum group $\text{\rm SU}_q(2)$,
{\it Internat.\ Math.\ Res.\ Notices} {\bf 2000}, 35--56.

\medskip\noindent
[U5] Y.\ Ueda, Fullness, Connes' $\chi$-groups, and ultra-products 
of amalgamated free products over Cartan subalgebras, 
{\it Trans.\ Amer.\ Math.\ Soc.} {bf 355} (2003), 349--371. 

\medskip\noindent
[VDN] D.\-V.\ Voiculescu, K.\-J.\ Dykema and A.\ Nica, Free Random
Variables, CRM Monograph Series, Vol.\ 1, Amer.\ Math.\ Soc.\ Providence,
RI, 1992.

\enddocument